\documentclass[
  a4paper, 
  reqno, 
  oneside, 
  11pt
]{amsart}

\usepackage[dvipsnames]{xcolor} 
\colorlet{cite}{BlueViolet}
\usepackage{tikz}  
\usepackage{rotating}
\usetikzlibrary{arrows,positioning}
\tikzset{ 
  baseline=-2.3pt,
  text height=1.5ex, text depth=0.25ex,
  >=stealth,
  node distance=2cm,
  mid/.style={fill=white,inner sep=2.5pt},
}

\usepackage{lmodern}
\usepackage{tikz-cd}
\usepackage{ amssymb, amsfonts}

\usepackage{graphicx,caption,subcaption}
\usepackage{fourier}
\usepackage{braket}  
\usepackage{microtype}

\usepackage[dvipsnames]{xcolor}
\usepackage[%
  bookmarks=true,			
  unicode=true,			
  pdftitle={mirrors},		%
  pdfauthor={Elizabeth Gasparim}{},	%
  pdfkeywords={}{adjoint orbit}{compactification}{Lefschetz fibration},	
  colorlinks=true,		
  linkcolor=Blue,			
  citecolor=cite,		
  filecolor=magenta,		
  urlcolor=RoyalBlue			
]{hyperref}				
\usepackage{cleveref}
\newcommand{\ce}{\mathrel{\mathop:}=}

\theoremstyle{plain}	
\newtheorem{theorem}{Theorem} 

\newtheorem*{theorem*}{Theorem}
\newtheorem*{lemma*}{Lemma}

\newtheorem*{corollary*}{Corollary}

\newtheorem*{definition*}{Definition}
\newtheorem*{conjecture*}{Conjecture}

\newtheorem*{example*}{Example}

\DeclareMathOperator{\Hom}{Hom}

\DeclareMathOperator{\im}{im}

\begin{document}
\author{Elizabeth Gasparim }
\address{E. Gasparim- Depto. Matem\'aticas, Univ. Cat\'olica del Norte, Antofagasta, Chile. etgasparim@gmail.com}
\title{Intrinsic mirrors for minimal adjoint orbits (ICM G\&T)}

\maketitle

This text  written for the ZAG volume\footnote{I thank Ivan Cheltsov for inviting me to contribute to this volume.} summarises  the Short Communication I presented at the 
Geometry and Topology Session of the  2022 International Congress of Mathematicians
 which took place in Copenhagen.

\section{Geometric Mirror Symmetry}

The original version of the Mirror Symmetry conjecture included the statement  that for each Calabi--Yau 
variety $X$ there exists  a dual Calabi--Yau variety $X^\vee$ such that if the diamond of $X$ is given below, 
then the diamond of $X^\vee$ is obtained from the one of $X$ by reflection on the 45 degree line, 
that is, by exchanging red and black, by flipping over the diagonal line passing thought the symbol for 
Serre duality. 
{\small 
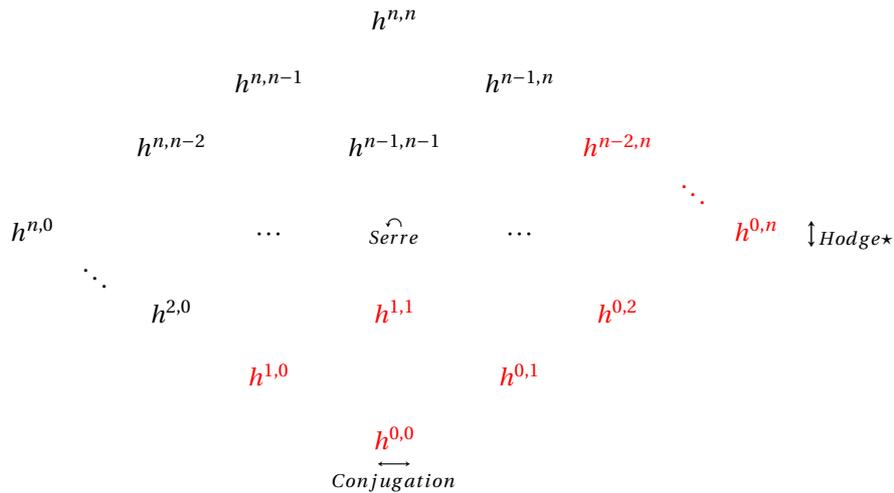
\begin{displaymath}
\begin{array}{cccccccccc}
 & & & & h^{n,n} & & & & & \\
 & & & & & & & & \\
 & & & h^{n,n-1} & & h^{n-1,n} & & & & \\
 & & & & & & & & & \\
 & & h^{n,n-2} & & h^{n-1,n-1} & & \textcolor{red}{h^{n-2,n}} & & & \\
 & & & &  & & & \textcolor{red}{\ddots} & & \\
 h^{n,0} & & & {\cdots} & {\stackrel{\curvearrowleft}{_{Serre}}}{} & {\cdots} & & & \textcolor{red}{ h^{0,n}} & 
 \updownarrow_{Hodge \star}\\
 & \ddots & & & {} & & & & & \\
 & & h^{2,0} & & \textcolor{red}{h^{1,1}} & & \textcolor{red}{h^{0,2}} & & & \\
 & & & & & & & & & \\
 & & & \textcolor{red}{h^{1,0}} & & \textcolor{red}{h^{0,1}} & & & & \\
 & & & & & & & & & \\
 & & & & \textcolor{red}{h^{0,0}} & & & & &\\
 & & & & \stackrel{\longleftrightarrow}{_{Conjugation}} & & & & &
\end{array}
\end{displaymath}
}

If a mirror $X^\vee$ failed to exist for $X$, then the variety $X$ would be said to behave like a {\it vampire}.
The validity of the conjecture would then imply that vampires do not exist. 


\section{Homological Mirror Symmetry Conjecture  I}

Maxim Kontsevitch presented the Homological Mirror Symmetry conjecture    at the at 1994 International Congress of Mathematicians:\\
 
\noindent{\bf HMS1}: For every  (subvariety of) a projective variety $X$  there exists a Landau--Ginzburg model $(Y,f)$ 
satisfying the categorical equivalence 
\textcolor{RoyalBlue}{$$Fuk(Y,f) \equiv D^b(Coh \, X)\text{.}$$}

\section{Symplectic Lefschetz fibrations}
Here, 
a { Landau--Ginzburg model} is a   manifold $Y$ together with a complex function 
$f\colon Y \rightarrow  \mathbb C$  (or $\mathbb P^1$)
called the superpotential. Typical examples are given by 
symplectic families with 1 dimensional parameter  space and only Morse type singularities, 
known as {\it symplectic Lefschetz fibrations}. Many examples of SLFs in dimension 4 were given by Donaldson
in his description of symplectic 4-manifolds,
but examples in higher dimensions were lacking. To provide  
examples of SLFs in higher dimensions we carried out a very general construction using methods from 
Lie theory \cite{GGS1,GGS2} which we now recall.

Let $G$ be a complex semisimple Lie group with Lie algebra $\mathfrak{g}$  and  Cartan 
subalgebra $\mathfrak h$. 
Given the Hermitian form 
 $\mathcal{H}$ on $\mathfrak g$, define the symplectic form on $\mathfrak g$ by 
$
\omega(X,Y) = \im \mathcal{H}\left( X_1,X_2\right).$
 For $H_0 \in \mathfrak h $ we consider the  adjoint orbit:
$
\mathcal{O}\left( H_0\right) =\mathrm{Ad}(G)\cdot H_0=
\{\mathrm{Ad}(g)\cdot H_0\in \mathfrak{g}:g\in G\}, $
together with the symplectic form $\omega$.

\begin{theorem*}[Gasparim, Grama, San Martin]
Given $H_{0}\in \mathfrak{h}$ and $H\in \mathfrak{h}_{%
\mathbb{R}}$ with $H$ a regular  element. The
\textquotedblleft height function\textquotedblright\ $f_{H}:(\mathcal{O}%
\left( H_{0}\right), \omega) \rightarrow \mathbb{C}$ defined by 

\centerline{\textcolor{RoyalBlue}{\boxed{$$
f_{H}\left( x\right) =\langle H,x\rangle \qquad x\in \mathcal{O}\left(
H_{0}\right) 
$$}}
}

\noindent has a finite number 
of isolated singularities and defines a  symplectic Lefschetz fibration.
\end{theorem*}

My favourite details of the proof:
$x$ is  a critical point for $f_{H}$ if and only if $x\in \mathcal{O}\left(
H_{0}\right) \cap \mathfrak{h}=\mathcal{W}\cdot H_{0}$, where $\mathcal{W}$ 
is the Weyl group. \\

For the purpose of calculating the Fukaya category $Fuk(Y,f)$ of vanishing cycles that appear in MHS1, 
it is important to find out whether   these SLFs allow for a rich variety of 
middle homology of its regular fibres. To understand this, we discuss the comparison 
between our semisimple adjoint orbits   and cotangent bundles of flag manifolds. 

Let $K \subset G$ compact, and $\mathbb F_0:= \mathrm{Ad}\left(
K\right) \cdot H_{0}$ be the  flag manifold.

\begin{theorem*}[Gasparim, Grama, San Martin]

There exists a $\mathbb C^\infty$  real  isomorphism $\iota \colon\mathrm{Ad}\left(
G\right) \cdot H_{0}\rightarrow T^{\ast }\mathbb{F}_{0}$ such that

\begin{enumerate}
\item $\iota $ is equivariant with respect to the action of 
$K$.

\item The pullback  of the canonical symplectic form on 
$T^{\ast }\mathbb{F}_{0}$ by $\iota$ is the  (real) Kirillov--Kostant--Souriau 
form on the  orbit.
\end{enumerate}
\end{theorem*}

As a consequence, we are  able to describe the topology of a regular level.

\begin{corollary*}
The homology of a  regular level of the SLF defined by $f_H$  coincides  with the homology of \,
$\mathbb{F}_{0}\setminus \mathcal{W}\cdot H_{0}$.
In particular the middle Betti number is $k-1$ where $k$ is 
the number of singularities of $f_H$ (equal the number of elements in
$\mathcal W \cdot H_0$).
\end{corollary*}

\section{The Katzarkov--Kontsevich--Pantev conjecture}

Katzarkov, Kontsevich, and Pantev defined 3 new Hodge theoretical invariants, 
which apply to a Landau--Ginzburg model, that is, which take into consideration not 
only the Hodge theory of $Y$ but also the potential function $f$. They then conjectured 
that these new 3 invariants coincide. Although the conjecture was proved to be false 
by Lunts and Przyjalkowski \cite{LP},  
Cheltsov and Przyjalkowski proved the KKP  conjecture for Fano threefolds \cite{CP}, and
we proved  it for some of our SLFs \cite{BGRS} (its validity for  all of our SLFs is as yet unknown).

Let  $\mathrm{LG}(n)$ denotes the Landau--Ginzburg model defined over 
the minimal semisimple adjoint orbit of $\mathfrak{sl}(n,\mathbb C)$, that is, the one
diffeomorphic with $T^*\mathbb P^n$.

\begin{theorem*}[Ballico, Gasparim, Rubilar, San Martin]
 The Landau--Ginzburg model $\mathrm{LG}(n)$ satisfies the KKP conjecture.
\end{theorem*}

Therefore, from the Hodge theoretical viewpoint, our LG models are behaving well. 
We may then proceed to the exciting task of computing the Fukaya category.

\section{The Fukaya category}

$Fuk(Y,f)$ is a  category whose objects are {\it Lagrangian thimbles}
associated to the vanishing cycles of the SLF.

\begin{definition*}The 
category of vanishing cycles $Lag(Y,f)$ is an 
$A_{\infty}$-category which objects $L_0,\dots, L_r$ corresponding to the 
thimbles.
The morphisms between the objects are given by 
$$\Hom(L_i,L_j)= \left\{\begin{array}{lll}
CF^*(L_i,L_j;R) = R^{[L_i\cap L_j]} & if & i<j \cr
R\cdot id & if & i=j \cr
0 & if & i>j
\end{array}\right. $$

The differential $m_1$, composition $m_2$ and higher order products $M_k$ 
are defined in terms of Lagrangian Floer homology.

\end{definition*}

\begin{example*}{ $\mathrm{LG}(2)$ is the Landau--Ginzburg model formed by 
 the adjoint orbit $\mathcal O_2$ of $\mathfrak{sl}\left( 2,\mathbb{C}\right) $,} and potential $f_H$
with the choices:
\[
H=H_{0}=\left( 
\begin{array}{cc}
1 & 0 \\ 
0 & -1%
\end{array}%
\right) . 
\]%
\begin{itemize}
\item Hence $\mathcal{O}_2$ is the set of 
matrices in  $\mathfrak{sl}\left( 2,\mathbb{C}\right) $ with eigenvalues $\pm 1$. 
\item  $\mathcal{O}_2$ forms a submanifold  of $\mathfrak{sl}\left( 2,\mathbb{C}\right) $ of real dimension $4$.
\item In this  case 
the potential $f_{H}=:\mathcal{O}_2\rightarrow \mathbb C$ has two singularities: $\pm H_0$. 
\end{itemize}

\begin{lemma*}
\label{teosl2} The Fukaya--Seidel category  \textcolor{ProcessBlue}{$Fuk(\mathrm{LG}(2))$}
is generated by two Lagrangians $L_0$ and $L_1$ with morphisms: 
\begin{align}
\label{morphismsfukaya}
\Hom (L_i, L_j) 
\simeq
\begin{cases}
\mathbb Z \oplus \mathbb Z [-1] & i < j \\
\mathbb Z                       & i = j \\
0                               & i > j
\end{cases}
\end{align}
where we think of $\mathbb Z$ as a complex concentrated in degree $0$ and $\mathbb Z [-1]$ as its shift, concentrated in degree $1$, and the products $m_k$ all vanish except for $m_2(\,\cdot\,, \mbox{id})$ and $m_2({\mbox{id},}\,\cdot\,)$.
\end{lemma*}

\begin{theorem*}[Ballico,  Barmeier,  Gasparim, Grama, San Martin] HMS1 fails for $\mathrm{LG}(2)$:

\begin{itemize}
\item $\mathrm{LG}(2)$ has no projective mirrors.

\item $\overline{\mathrm{LG}(2)}$ has no projective mirrors. 
\end{itemize}
\end{theorem*}
This result proved in \cite{BBGGS} means that
for any (subvariety of a) projective variety $X$  we have 
 $$\textcolor{ProcessBlue}{Fuk(\mathrm{LG}(2))} \not\equiv D^bCoh(X).   $$

\end{example*}

\section{Homological Mirror Symmetry Conjecture  II}

\noindent{\bf HMS2}: For every LG model $(Y,f)$ there exists an LG model  $(X,g)$ 
such that 
\textcolor{RoyalBlue}{$$Fuk(Y,f) \equiv D_{Sg}(X,g)$$}
where $D_{Sg}(X,g)$ denotes the Orlov category of singularities.

\begin{definition*}The { Orlov category of  singularity  } of $(X,g)$ is 

$$D_{Sg}(X,g)\ce\bigoplus_i \frac{D^bCoh (X_i)}{\mathfrak{Perf}(X_i)}$$
where $X_i$ are the critical fibers of $g$ and perfect complexes are  (quasi-isomorphic to)
those of the form:
$$E_1 \rightarrow E_2 \rightarrow \cdots \rightarrow E_n $$
where $E_i$ are locally free, see \cite{O1,O2}.
\end{definition*}

\begin{definition*}The {Intrinsic Mirror Symmetry } program of  Gross
and Siebert \cite{GS} provides a  recipe to find the mirror:
\begin{itemize}
\item Fit the LG model inside a log Calabi--Yau pair. 
\item Compute the intersection complex for the dual pair. 
\item Construct Theta functions.
\item Compute punctured Gromov--Witten invariants.
\end{itemize}
\end{definition*}

  Using the Gross--Siebert  recipe, in  \cite{G} I obtained  a Landau--Ginzburg model 
   mirror to $\mathrm{LG}(2)$ as
$\mathrm{LG}^\vee\!(2)\ce (X_2,g=y)$ where $X_2 \subset \mathbb C \times \mathbb C^*\times \mathbb P^1$
is given by the equation:
$$
uy=v(x+1+1/x).
$$

 \begin{theorem*}[HMS2 for $\mathrm{LG}(2)$]There is an equivalence of categories:
 $$\textcolor{ProcessBlue}{Fuk(\mathrm{LG}(2)) \equiv \textcolor{ProcessBlue}{D_{Sg}\mathrm{LG}^\vee\!(2)} }.$$
 \end{theorem*}

\end{document}